# A General Conditional Large Deviation Principle

Brian R. La Cour[1] · William C. Schieve[2]

**Abstract** Given a sequence of Borel probability measures on a Hausdorff space which satisfy a large deviation principle (LDP), we consider the corresponding sequence of measures formed by conditioning on a set $B$. If the large deviation rate function $I$ is good and effectively continuous, and the conditioning set has the property that (1) $\overline{B^\circ} = \overline{B}$ and (2) $I(x) < \infty$ for all $x \in \overline{B}$, then the sequence of conditional measures satisfies a LDP with the good, effectively continuous rate function $I_B$, where $I_B(x) = I(x) - \inf I(B)$ if $x \in \overline{B}$ and $I_B(x) = \infty$ otherwise.

## 1 Introduction

Conditional large deviation principles (LDPs) have played an important role in both mathematical statistics and statistical mechanics. In the latter field in particular, the asymptotic behavior of conditional distributions has been studied since the time of Boltzmann [1,5] and remains an important subject today in establishing the equivalence of microcanonical and canonical ensembles [6,14]. Early work by Lanford [9] and Ruelle [13] had anticipated many of the results regarding entropy and the thermodynamic limit which would later find form in the mathematical theory of large deviations. More recently, large deviation techniques have also been applied to nonequilibrium macroscopic time evolution [7,8].

Much of the past work on conditional LDPs has focused on distributions of sample means or, equivalently, empirical measures for mutually independent and identically distributed (IID) random variables. In this case, the conditioning set may be expressed as a constraint on

B  Brian R. La Cour
   blacour@arlut.utexas.edu

[1] Applied Research Laboratories, The University of Texas at Austin, P. O. Box 8029, Austin, TX 78713-8029, USA

[2] Department of Physics, The University of Texas at Austin, 1 University Station C1600, Austin, TX 78712, USA

the sample mean of IID random variables or, equivalently, as a constraint on the empirical distribution directly. Stroock and Zeitouni [14], for example, have developed a theorem on Gibbs conditioning, later revised in Dembo and Zeitouni [3], which establishes convergence in probability for empirical measures conditioned on the sample mean of IID random variables. Earlier, van Campenhout and Cover [15], using results from Zabell and Lanford, showed explicitly that the marginals of a distribution conditioned on a particular value of the sample mean converge to the canonical form predicted by the maximum entropy principle.

General conditional limit theorems have also been studied by Lewis et al. [10]. Casting the work of Ruelle and Lanford into the language of large deviation theory, they show that for a suitable sequence of conditioning sets, the corresponding conditional measures converge asymptotically to a canonical or "tilted" (in the Varadhan sense) measure with respect to the Kullback–Leibler information gain (Theorem 5.1). If the conditioning sets are assumed to be convex, then an even stronger conclusion follows in which the conditional equilibrium points may be determined from the subgradient of the free energy (Theorem 6.1). Their results fall short of giving an LDP for the sequence of conditional measures, however, and it is this problem which we address here.

In this paper we develop a general conditional LDP in terms of an assumed LDP and a given conditioning set. Suppose $\{P_n\}$ is a sequence of Borel probability measures on a Hausdorff topological space $(X, \mathcal{T})$ which satisfies a LDP with a good rate function $I$. Given a set $B \subseteq X$ for which $\inf I(B^\circ) < \infty$, the large deviation lower bound implies that $P_n(B) > 0$ for all $n$ sufficiently large. Without loss of generality, we may suppose that $P_n(B) > 0$ for all $n$, so the sequence $\{P_n(\cdot|B)\}$ of conditional probability measures is well defined. Since the minimum of the unconditioned rate function, $I$, gives the asymptotically most likely value (or values) in $X$, we may anticipate that, under conditioning, the asymptotically most likely values in $B$ (or, more precisely, in $\overline{B}$) will be those values that minimize $I$ over $\overline{B}$. If this sequence of conditional measures is to have an LDP, however, then this minimum must be zero, and this suggests that the conditional rate function, denoted $I_B$, is given by $I_B(x) = I(x) - \inf I(\overline{B})$ for all $x \in \overline{B}$. For consistency, $I_B(x)$ may be defined to be $\infty$ outside of $\overline{B}$. Since $\inf I(\overline{B}) \leq \inf I(B^\circ)$ and $\inf I(B^\circ) < \infty$ by assumption, we see that $I_B$ is indeed well defined.

The possible discrepancy between $\inf I(B^\circ)$ and $\inf I(\overline{B})$ is an inconvenience which may be eliminated if reasonable regularity conditions are placed on $B$ and $I$. Equality of $\inf I(B^\circ)$ and $\inf I(\overline{B})$ of course holds for $I$-continuity sets, and this, in turn, holds if $\overline{B^\circ} = \overline{B}$ and $I$ is everywhere finite and continuous. (See Theorem 2 below.) The former condition is clearly satisfied if $B$ is open or, e.g., a closed ball in a metric space. Also, if $B$ is a convex subset of a normed linear space, then $\overline{B^\circ} = \overline{B}$ holds if $B$ has a nonempty interior. This is a straightforward extension of Theorem 6.3 in Rockafellar [12]. The utility of convex conditioning sets had been recognized by Csiszár [2] in studying conditionally distributed empirical measures. The condition that $I$ be finite and continuous everywhere may be further relaxed by noting that the latter condition is relevant only on the effective domain, $D_I = \{x \in X : I(x) < \infty\}$, of $I$. We shall say that a rate function is *effectively continuous* if it is continuous relative to its effective domain. Convex rate functions on a Banach space, for example, are effectively continuous, since $D_I$ is convex whenever $I$ is convex (see, e.g., Roberts and Varberg [11], p. 112). With these considerations in mind, the main result may now be stated as follows:

**Theorem 1** (Conditional LDP) *Suppose $\{P_n\}$ satisfies an LDP with a good rate function $I$ on a Hausdorff space $X$. Let $B$ be a given Borel set for which $\overline{B^\circ} = \overline{B}$ and $\emptyset \subset B^\circ \subseteq D_I$. If $I$ is continuous on $B^\circ$, then the sequence $\{P_n(\cdot|B)\}$ of conditional probability measures satisfies an LDP with the good rate function*

$$I_B(x) = \begin{cases} I(x) - \inf I(B), & \text{if } x \in \overline{B}, \\ \infty & \text{otherwise.} \end{cases} \quad (1)$$

This result is proved in Sect. 3. The restriction to Hausdorff spaces is needed to obtain a good conditional rate function.

## 2 Auxilliary Rate Function Theorems

Before proving the conditional LDP theorem, a few general results regarding large deviation rate functions are needed. We begin with the following theorem regarding continuous rate functions. Throughout this section it is assumed that $(X, \mathcal{T})$ is a Hausdorff topological space.

**Theorem 2** *Let $I$ be a good rate function which is continuous on an open set containing $A$. If either $A = \emptyset$ or $\inf I(A) < \infty$, then $\inf I(A) = \inf I(\overline{A})$.*

*Proof* If $A = \emptyset$ then $I$ is trivially continuous on $A$ and $\inf I(A) = \infty = \inf I(\overline{A})$. Now suppose $A \neq \emptyset$. Since $I$ is a good rate function and $\overline{A}$ is closed, there exists at least one $x_A \in \overline{A}$ such that $I(x_A) = \inf I(\overline{A})$. Since $A \subseteq \overline{A}$ and is nonempty, $I(x_A) \leq \inf I(A) < \infty$. Suppose $I(x_A) < \inf I(A)$. Now let $V = (-\infty, \inf I(A))$ and note that, since $I$ is continuous and $I(x_A) \in V$ by assumption, there exists a neighborhood $U$ of $x_A$ such that $I(U) \subseteq V$. Since $x_A \in \overline{A}$ and $U$ is a neighborhood of $x_A$, there exists an $x_A' \in U \cap A$. Since $x_A' \in U$, $I(x_A') \in V$ and hence $I(x_A') < \inf I(A)$; however, since $x_A' \in A$, $I(x_A') \geq \inf I(A)$. We thus arrive at a contradiction and conclude that $I(x_A) = \inf I(\overline{A}) = \inf I(A)$. □

**Lemma 1** *For any subsets $A$ and $B$ of a topological space, $A^\circ \cap \overline{B} \subseteq \overline{A^\circ \cap B}$.*

*Proof* If $A^\circ \cap \overline{B} = \emptyset$ then we are done, so suppose there exists an $x \in A^\circ \cap \overline{B}$. We will have $x \in \overline{A^\circ \cap B}$ if and only if for any neighborhood $U$ of $x$ we have that $U \cap (A^\circ \cap B) \neq \emptyset$. Now, given $U$, $U \cap A^\circ$ is also a neighborhood of $x$, then, since $x \in \overline{B}$, it follows that $(U \cap A^\circ) \cap B \neq \emptyset$. □

**Corollary 1** *Let $I$ be a good rate function which is continuous on $A^\circ \cap B^\circ$. If $\overline{B^\circ} = \overline{B}$ and either $A^\circ \cap B^\circ = \emptyset$ or $\inf I(A^\circ \cap B^\circ) < \infty$, then $\inf I(A^\circ \cap B^\circ) = \inf I(A^\circ \cap \overline{B})$.*

*Proof* Since $\overline{B} = \overline{B^\circ}$ and $A^\circ \cap \overline{B^\circ} \subseteq \overline{A^\circ \cap B^\circ}$ by Lemma 1, we have

$$\inf I(A^\circ \cap B^\circ) \geq \inf I(A^\circ \cap \overline{B}) = \inf I(A^\circ \cap \overline{B^\circ}) \geq \inf I(\overline{A^\circ \cap B^\circ}) = \inf I(A^\circ \cap B^\circ),$$

where the last equality follows from Theorem 2. □

## 3 Proof of Conditional Large Deviation Principle

Observe that the large deviation bounds imply that, for any $\varepsilon > 0$ and all $n$ sufficiently large,

$$a_n^{-1} \log P_n(A) < -(1-\varepsilon) \inf I(\overline{A}) \quad \text{for } 0 < \inf I(\overline{A}) < \infty, \quad (2)$$
$$a_n^{-1} \log P_n(A) > -(1+\varepsilon) \inf I(A^\circ) \quad \text{for } 0 < \inf I(A^\circ) < \infty, \quad (3)$$

where $\{a_n\}$ is an unbounded sequence of positive scale factors. Similarly, $\inf I(\overline{A}) = \infty$ implies $P_n(A) = 0$ for all $n$ sufficiently large, while $\inf I(A^\circ) = 0$ implies $a_n^{-1} \log P_n(A) > -\varepsilon$ for $\varepsilon > 0$ and all $n$ sufficiently large. With these observations, we are now ready to prove Theorem 1.

*Proof* Since $I$ is a good rate function which is continuous on $B°$ and $B° \subseteq D_I$, $\infty >$ inf $I(B°) = $ inf $I(\overline{B°})$, by Theorem 2. Furthermore, since $\overline{B°} = \overline{B}$, inf $I(\overline{B°}) = $ inf $I(\overline{B}) = $ inf $I(B)$. If inf $I(B) > 0$, then

$$-(1+\varepsilon) \inf I(B) < a_n^{-1} \log P_n(B) < -(1-\varepsilon) \inf I(B)$$

for $\varepsilon > 0$ and all $n$ sufficiently large, while inf $I(B) = 0$ implies

$$-\varepsilon < a_n^{-1} \log P_n(B) \leq -(1-\varepsilon) \inf I(B) = 0.$$

We begin with the large deviation upper bound. First assume $0 < \inf I(\overline{A \cap B}) < \infty$ and inf $I(B) > 0$. For a given $\varepsilon > 0$ we have that for all $n$ sufficiently large

$$\begin{aligned} a_n^{-1} \log P_n(A|B) &= a_n^{-1} \log P_n(A \cap B) - a_n^{-1} \log P_n(B) \\ &< -(1-\varepsilon) \inf I(\overline{A \cap B}) + (1+\varepsilon) \inf I(B°) \\ &\leq -\left[\inf I(\overline{A} \cap \overline{B}) - \inf I(B)\right] + \varepsilon \left[\inf I(\overline{A \cap B}) + \inf I(B)\right] \\ &= -\inf I_B(\overline{A}) + \varepsilon \left[\inf I(\overline{A \cap B}) + \inf I(B)\right] \end{aligned}$$

As the second term is positive and may be made arbitrarily small, we conclude

$$\limsup_{n \to \infty} a_n^{-1} \log P_n(A|B) \leq -\inf I_B(\overline{A}).$$

If $0 < \inf I(\overline{A \cap B}) < \infty$ yet inf $I(B) = 0$, then

$$\begin{aligned} a_n^{-1} \log P_n(A|B) &< -(1-\varepsilon) \inf I(\overline{A \cap B}) + \varepsilon \\ &\leq -\inf I_B(\overline{A}) + \varepsilon \left[\inf I(\overline{A \cap B}) + 1\right] \end{aligned}$$

and the upper bound is again found to hold.

If inf $I(\overline{A \cap B}) = 0$, then inf $I(B) = $ inf $I(\overline{B}) \leq $ inf $I(\overline{A} \cap \overline{B}) \leq $ inf $I(\overline{A \cap B}) = 0$ and inf $I_B(\overline{A}) = $ inf $I(\overline{A} \cap \overline{B}) - $ inf $I(B) \leq $ inf $I(\overline{A \cap B}) - $ inf $I(B) = 0$. Since $a_n^{-1} \log P_n(A|B) \leq 0 = -\inf I_B(\overline{A})$, the upper bound is clearly satisfied.

If inf $I(\overline{A \cap B}) = \infty$, then $P_n(A|B) = 0$ and $a_n^{-1} \log P_n(A|B) = -\infty$ for all $n$ sufficiently large. Thus, $\limsup_{n \to \infty} a_n^{-1} \log P_n(A|B) = -\infty \leq -\inf I_B(\overline{A})$.

For the large deviation lower bound, suppose $0 < \inf I((A \cap B)°) < \infty$ and note that for all $n$ sufficiently large,

$$\begin{aligned} a_n^{-1} \log P_n(A|B) &> -(1+\varepsilon) \inf I((A \cap B)°) + (1-\varepsilon) \inf I(B) \\ &= -\left[\inf I(A° \cap B°) - \inf I(B)\right] - \varepsilon \left[\inf I(A° \cap B°) + \inf I(B)\right] \\ &= -\inf I_B(A°) - \varepsilon \left[\inf I(A° \cap B°) + \inf I(B)\right], \end{aligned}$$

where Corollary 1 has been used in the last equality. The second term is negative and may be made arbitrarily small, so we conclude

$$\liminf_{n \to \infty} a_n^{-1} \log P_n(A|B) \geq -\inf I_B(A°).$$

If inf $I(A° \cap B°) = 0$, then inf $I(B) = $ inf $I(\overline{B}) \leq $ inf $I(A° \cap \overline{B}) = $ inf $I(A° \cap B°) = 0$ and inf $I_B(A°) = $ inf $I(A° \cap \overline{B}) - $ inf $I(B) = $ inf $I(A° \cap B°) - $ inf $I(B) = 0$. However, for any given $\varepsilon > 0$ and all $n$ sufficiently large,

$$a_n^{-1} \log P_n(A|B) > -\varepsilon + (1-\varepsilon) \inf I(B) = -\varepsilon = -\inf I_B(A°) - \varepsilon.$$

Thus, the lower bound is satisfied in this case.

Finally, suppose $\inf I(A° \cap B°) = \infty$. Since $B° \subseteq D_I$, this implies $A° \cap B° = \emptyset$. By Corollary 1, $\inf I_B(A°) = \inf I(A° \cap \overline{B}) - \inf I(B) = \inf I(A° \cap B°) - \inf I(B) = \infty$. But since $a_n^{-1} \log P_n(A|B) \geq -\infty = -\inf I_B(A° \cap B°)$, the lower bound is clearly satisfied in this case as well.

To complete the proof, we must show that $I_B$ is a good, continuous rate function relative to $\overline{B}$. Effective continuity of $I_B$ follows from that of $I$. To show that $I_B$ is a good rate function, consider any $\alpha < \infty$ and note that

$$\{x \in X : I_B(x) \leq \alpha\} = \{x \in \overline{B} : I(x) - \inf I(B) \leq \alpha\}$$
$$= \{x \in X : I(x) \leq \alpha + \inf I(B)\} \cap \overline{B}.$$

We have already established that $\inf I(B) < \infty$. As $X$ is a Hausdorff space and $I$ is a good rate function, the above intersection is compact, thus establishing that $I_B$ is a good rate function. □

## 4 Application to Joint Random Vectors

Let $\{(\Omega_n, \mathcal{F}_n, P_n)\}_{n \in \mathbb{N}}$ be a sequence of probability spaces and let $\{(X_n, Y_n)\}_{n \in \mathbb{N}}$ be a sequence of Borel-measurable random vectors on $\Omega_n$ taking values in $\mathbb{R}^d \times \mathbb{R}^{d'}$. Suppose we are interested in the asymptotic behavior of $Y_n$ when $X_n$ is conditioned on a value $x_0 \in \mathbb{R}^d$. Rather than condition on $X_n = x_0$ explicitly, we shall instead consider the joint distribution of $(X_n, Y_n)$ and construct a conditioning set for which $X_n$ converges to $x_0$ in probability. Assuming an LDP for the joint distribution and using the conditional LDP theorem (Theorem 1), we will determine the value $y_0$ corresponding to $x_0$ to which $Y_n$ converges in probability under this conditioning.

By the Gärtner–Ellis Theorem [3], the joint distribution of $(X_n, Y_n)$ will satisfy an LDP if the free energy, $\Psi : \mathbb{R}^d \times \mathbb{R}^{d'} \to (-\infty, \infty]$, given by

$$\Psi(\lambda_1, \lambda_2) = \lim_{n \to \infty} \frac{1}{a_n} \log \int_{\Omega_n} e^{a_n[\lambda_1 \cdot X_n(\omega) + \lambda_2 \cdot Y_n(\omega)]} dP_n(\omega) \quad (4)$$

is well defined and everywhere finite and differentiable. Assuming this to be the case, the rate function, $I$, is given by the Legendre-Fenchel transform of $\Psi$, i.e.,

$$I(x, y) = \lambda_x \cdot x + \lambda_y \cdot y - \Psi(\lambda_x, \lambda_y), \quad (5)$$

where $\lambda_x$ and $\lambda_y$ are such that $x = \nabla_1 \Psi(\lambda_x, \lambda_y)$ and $y = \nabla_2 \Psi(\lambda_x, \lambda_y)$. (The mapping $(\lambda_x, \lambda_y) \mapsto (x, y)$ is invertible if the Jacobian of $(\nabla_1 \Psi, \nabla_2 \Psi)$ exists everywhere and vanishes nowhere, by the inverse function theorem; we shall assume this is indeed the case.) It follows that $I$ is a good, essentially strictly convex (hence, effectively continuous) rate function. (See Ellis [4], Theorem VII.2.1.) As such, there is a unique point $(x_*, y_*)$ for which the rate function is zero and to which $(X_n, Y_n)$ converges in probability. In terms of the free energy, note that $x_* = \nabla_1 \Psi(0, 0)$ and $y_* = \nabla_2 \Psi(0, 0)$. The effective domain of $I$ will be denoted, as usual, by $D_I$.

To condition on $x_0$, we consider a conditional LDP with a conditioning set $B$ chosen so that $I$ has its infimum at a unique point $(x_0, y_0)$ for some $y_0$. Not all choices of $x_0$ will allow for a suitable conditioning set, as boundary points may be problematic. One way to address this problem is to consider the LDP for $X_n$ alone. By the contraction principle, $X_n$ satisfies an LDP with rate function $I_X(\cdot) = I(\cdot, y_*)$ and corresponding free energy $\Psi_X(\cdot) = \Psi(\cdot, 0)$. If we choose $x_0 \in \nabla \Psi_X(\mathbb{R}^d)$, then clearly there exists a $\lambda_0 \in \mathbb{R}^d$ such that $x_0 = \nabla \Psi_X(\lambda_0) = \nabla_1 \Psi(\lambda_0, 0)$. As we have assumed invertibility, $\lambda_0$ is in fact uniquely determined by $x_0$. Now choose

$$B = \{(x, y) \in D_I : \lambda_0 \cdot (x - x_0) \geq 0\}. \tag{6}$$

Using the value of $\lambda_0$ determined by $x_0$, define $y_0 = \nabla_2 \Psi(\lambda_0, 0)$ and note that, since $I(x_0, y_0) = \lambda_0 \cdot x_0 + 0 \cdot y_0 - \Psi(\lambda_0, 0) < \infty$, $(x_0, y_0) \in D_I$. From its definition, $B$ is the intersection of the convex set $D_I$ with the affine half-space demarcated by the hyperplane containing the point $(x_0, y_0)$ and having a normal vector proportional to $(\lambda_0, 0)$ directed towards its interior; thus, $B$ is also convex. We shall now verify that the conditions of Theorem 1 do indeed hold.

We have already established that $I$ is a good, effectively continuous rate function and its domain is clearly Hausdorff, so it remains to verify the required conditions on $B$. Clearly $B$ is a convex set, and $B^\circ \subseteq B \subseteq D_I$. Due to the choice of $x_0$ and the assumed continuity of $(\nabla_1 \Psi, \nabla_2 \Psi)$, $B$ also has a nonempty interior, so the fact that it is convex implies $\overline{B^\circ} = \overline{B}$. This establishes the conditional LDP. It remains, then, to determine the corresponding rate function, i.e., to compute $\inf I(B)$.

For $(x, y) \in B$ we have that $\lambda_0 \cdot x \geq \lambda_0 \cdot x_0$, so

$$\begin{aligned} I(x, y) &= \lambda_x \cdot x + \lambda_y \cdot y - \Psi(\lambda_x, \lambda_y) \\ &\geq \left[\lambda_x \cdot x + \lambda_y \cdot y - \Psi(\lambda_x, \lambda_y) - \lambda_0 \cdot x + \Psi(\lambda_0, 0)\right] + I(x_0, y_0), \end{aligned}$$

since $\lambda_0 \cdot x - \Psi(\lambda_0, 0) \geq \lambda_0 \cdot x_0 - \Psi(\lambda_0, 0) = I(x_0, y_0)$. The expression in brackets is nonnegative, since

$$\begin{aligned} \lambda_x \cdot x + \lambda_y \cdot y - \Psi(\lambda_x, \lambda_y) &= \sup_{\lambda_1, \lambda_2} \left[\lambda_1 \cdot x + \lambda_2 \cdot y - \Psi(\lambda_1, \lambda_2)\right] \\ &\geq \lambda_0 \cdot x + 0 \cdot y - \Psi(\lambda_0, 0) \end{aligned}$$

with equality if and only if $(\lambda_x, \lambda_y) = (\lambda_0, 0)$. Thus, $I(x, y) \geq I(x_0, y_0)$ for all $(x, y) \in B$, and, since $(x_0, y_0) \in B$, we conclude that

$$\inf I(B) = I(x_0, y_0) = \lambda_0 \cdot x_0 - \Psi(\lambda_0, 0) \tag{7}$$

with $y_0 = \nabla_2 \Psi(\lambda_0, 0)$ and $\lambda_0$ given by $x_0 = \nabla_1 \Psi(\lambda_0, 0)$. This gives the desired conditional rate function, from which it follows that $(X_n, Y_n)$ converges in probability to $(x_0, y_0)$. Note that, since $B \subseteq D_I$ is convex and $I$ is essentially strictly convex, $(x_0, y_0)$ is the unique point in $B$ at which $I$ attains the minimum value $\inf I(B)$.

Note that this result continues to hold if the condition $\lambda_0 \cdot (x - x_0) \geq 0$ is replaced with $0 \leq \lambda_0 \cdot (x - x_0) < \delta$, where $\delta > 0$ is arbitrary. (In statistical mechanics, this corresponds to a microcanonical distribution with a "thickened" energy shell.) Consequently, while we do not condition on $X_n = x_0$ precisely, we may restrict $X_n$ to be arbitrarily close to $x_0$. The asymptotic value for $Y_n$, i.e., $y_0$, may be written in a more familiar form by evaluating $\nabla_2 \Psi(\lambda_0, 0)$ explicitly. Since $\Psi$ is convex and we have assumed it to be finite and differentiable, it follows from Theorem 25.7 of Rockafellar [12] that the convergence of the gradients is uniform; hence,

$$y_0 = \nabla_2 \Psi(\lambda_0, 0) = \lim_{n \to \infty} \frac{\int_{\Omega_n} Y_n(\omega) e^{a_n \lambda_0 \cdot X_n(\omega)} dP_n(\omega)}{\int_{\Omega_n} e^{a_n \lambda_0 \cdot X_n(\omega')} dP_n(\omega')}, \tag{8}$$

which is the familiar canonical expectation.

Using the contraction principle one may obtain an explicit LDP for $Y_n|_{x_0}$, i.e., $Y_n$ conditioned on $\lambda_0 \cdot (X_n - x_0) \in [0, \delta)$. To do this, note that the projection map $(x, y) \mapsto y$ is continuous; hence, $Y_n|_{x_0}$ satisfies an LDP with rate function

$$\begin{aligned} I_{x_0}(y) &= \inf\left\{I_B(x', y') : x' \in \mathbb{R}^d, y' = y\right\} \\ &= \inf\left\{I(x', y) : x' \in \mathbb{R}^d, \lambda_0 \cdot (x' - x_0) \in [0, \delta)\right\} - \inf I(B) \\ &= I(x_0, y) - \lambda_0 \cdot x_0 + \Psi(\lambda_0, 0) \end{aligned} \quad (9)$$

for $y \in \mathbb{R}^{d'}$. The last equality follows from an argument similar to that used to determine $\inf I(B)$. From the properties of $I$ it follows that $I_{x_0}$ is a good, essentially strictly convex rate function. Thus, the corresponding free energy, call it $\Psi_{x_0}$, is given by

$$\begin{aligned} \Psi_{x_0}(\lambda) &= \sup_{y \in \mathbb{R}^{d'}} \left[\lambda \cdot y - I_{x_0}(y)\right] \\ &= \sup_{y \in \mathbb{R}^{d'}} \left[\lambda_0 \cdot x_0 + \lambda \cdot y - I(x_0, y) - \Psi(\lambda_0, 0)\right] \\ &= \Psi(\lambda_0, \lambda) - \Psi(\lambda_0, 0) \end{aligned} \quad (10)$$

for $\lambda \in \mathbb{R}^{d'}$. The properties of finiteness and differentiability for $\Psi_{x_0}$ follow from those of $\Psi$ (or, more specifically, from those of $\Psi(\lambda_0, \cdot)$), so the global minimum of $I_{x_0}$ is attained at $\nabla \Psi_{x_0}(0) = \nabla_2 \Psi(\lambda_0, 0) = y_0$, as expected. Note that $I_{x_0}$ may also be written directly in terms of $\Psi_{x_0}$ via the relation

$$I_{x_0}(y) = \lambda \cdot y - \Psi_{x_0}(\lambda). \quad (11)$$

**Acknowledgments** This work was supported by the Engineering Research Program of the Office of Basic Energy Sciences at the U.S. Department of Energy, Grant DE-FG03-94ER14465. This work was also supported in part by Applied Research Laboratories, The University of Texas at Austin, under an Independent Research grant. The authors would like to thank R. Ellis for kindly providing a preprint of reference [6].